\documentclass{article}

\usepackage{microtype}
\usepackage{graphicx}
\usepackage{subfigure}
\usepackage{booktabs} 

\usepackage{hyperref}
\usepackage{ulem}
\newcommand{\stkout}[1]{\ifmmode\text{\sout{\ensuremath{#1}}}\else\sout{#1}\fi}

\newcommand{\RNum}[1]{\uppercase\expandafter{\romannumeral #1\relax}}

\usepackage[accepted]{icml2024}

\usepackage{amsmath}
\usepackage{amssymb}
\usepackage{mathtools}
\usepackage{amsthm}

\usepackage[capitalize,noabbrev]{cleveref}

\theoremstyle{plain}
\newtheorem{theorem}{Theorem}[section]

\theoremstyle{definition}
\newtheorem{definition}[theorem]{Property}

\theoremstyle{remark}

\usepackage[textsize=tiny]{todonotes}

\icmltitlerunning{Structure-Preserving Operator Learning: Modeling the Collision Operator of Kinetic Equations}

\newcommand{\Stwo}{{\mathbf{S}^2}}

\newcommand{\norm}[1]{\left\lVert#1\right\rVert}
\newcommand{\abs}[1]{\left|#1\right|}
\newcommand{\set}[1]{\left\lbrace#1\right\rbrace}

\newcommand{\intD}{\;\mathrm{d}}
\newcommand{\rom}[1]{\uppercase\expandafter{\romannumeral #1\relax}}

\newcommand{\inner}[1]{\left< #1 \right>}

\newcommand{\fu}{\mathbf{u}}

\newcommand{\fv}{\mathbf{v}}

\newcommand{\fx}{\mathbf{x}}
\newcommand{\fm}{\mathbf{m}}

\newcommand{\balpha}{{\boldsymbol{\alpha}}}

\begin{document}

\twocolumn[
\icmltitle{Structure-Preserving Operator Learning: Modeling the Collision Operator of Kinetic Equations}



\icmlsetsymbol{equal}{*}

\begin{icmlauthorlist}
\icmlauthor{Jae Yong Lee}{equal,kias}
\icmlauthor{Steffen Schotth\"ofer}{equal,ornl}
\icmlauthor{Tianbai Xiao}{equal,imech,ucas}
\icmlauthor{Sebastian Krumscheid}{kit}
\icmlauthor{Martin Frank}{kit}
\end{icmlauthorlist}

\icmlaffiliation{kit}{Karlsruhe Institute of Technology, Scientific Computing Center, 76131 Karlsruhe, Germany}
\icmlaffiliation{ornl}{Oak Ridge National Laboratory, Computer Science and Mathematics Division, Oak Ridge, Tennessee, USA}
\icmlaffiliation{kias}{Center for Artificial Intelligence and Natural Sciences, Korea Institute for Advanced Study,
Seoul, Republic of Korea}
\icmlaffiliation{imech}{State Key Laboratory of High Temperature Gas Dynamics and Centre for Interdisciplinary Research in Fluids, Institute of Mechanics, Chinese Academy of Sciences, Beijing, China}
\icmlaffiliation{ucas}{School of Engineering Science, University of Chinese Academy of Sciences, Beijing, China}

\icmlcorrespondingauthor{Steffen Schotth\"ofer}{schotthofers@ornl.gov}

\icmlkeywords{Machine Learning, ICML}

\vskip 0.3in
]



\printAffiliationsAndNotice{\icmlEqualContribution} 

\begin{abstract}
This work explores the application of deep operator learning principles to a problem in statistical physics. Specifically, we consider the linear kinetic equation, consisting of a differential advection operator and an integral collision operator, which is a powerful yet expensive mathematical model for interacting particle systems with ample applications, e.g., in radiation transport. We investigate the capabilities of the Deep Operator network (DeepONet) approach to modelling the high dimensional collision operator of the linear kinetic equation. This integral operator has crucial analytical structures that a surrogate model, e.g., a DeepONet, needs to preserve to enable meaningful physical simulation. We propose several DeepONet modifications to encapsulate essential structural properties of this integral operator in a DeepONet model. To be precise, we adapt the architecture of the trunk-net so the DeepONet has the same collision invariants as the theoretical kinetic collision operator, thus preserving conserved quantities, e.g., mass, of the modeled many-particle system.
Further, we propose an entropy-inspired data-sampling method tailored to train the modified DeepONet surrogates without requiring an excessive expensive simulation-based data generation.\footnote{
This manuscript has been authored by UT-Battelle, LLC under Contract No. DE-AC05-00OR22725 with the U.S. Department of Energy. The United States Government retains and the publisher, by accepting the article for publication, acknowledges that the United States Government retains a non-exclusive, paid-up, irrevocable, world-wide license to publish or reproduce the published form of this manuscript, or allow others to do so, for United States Government purposes. The Department of Energy will provide public access to these results of federally sponsored research in accordance with the DOE Public Access Plan(\url{http://energy.gov/downloads/doe-public-access-plan}).}
\end{abstract}

\section{Introduction}\label{sec:intro}

\subsection{General introduction and problem setting}
Statistical physics uses a statistical ensemble approach to model the behavior of interacting particle systems.
As a direct consequence of the ensemble model, the Boltzmann-type kinetic equation provides the mean-field description by tracking the space-time evolution of the particle probability distribution in the physical-velocity phase space.
Due to the high dimensionality and nonlinearity, obtaining accurate solutions to the kinetic equation is intricate and computationally demanding, even with modern data streaming architecture.
Bringing in scientific machine learning techniques, especially operator learning, is expected to significantly improve the computational efficiency of the kinetic equation and promote its further application in, e.g., aerospace~\cite{BoltzmannApplications}, the semiconductor industry~\cite{Markowich1990SemiconductorE}, and radiation transport~\cite{NeutronTransport, Chahine1987FoundationsOR}.
The kinetic equation takes the form
\begin{equation}  \label{eq_boltzmann}
	\partial_t f+\mathbf{v} \cdot \nabla_{\fx} f +\sigma_{\textup{a}} = \sigma_{\textup{s}}Q(f) + p,
\end{equation}
where $f(t,\fx, \fv)$ denotes the kinetic density in physical space $\fx\in \mathbf X\subset\mathbb{R}^3$, and velocity space $\mathbf{v}\in \mathbb{R}^3$, at time instant $t$.
The left-hand side of the kinetic equation describes the effects of particle flight due to density gradient, and the right-hand side models the interaction with the background medium. The term  $p(\fx,\fv)$ models particle sources, $\sigma_{\textup{a}}(\fx)$ models absorption by the background material, and $ \sigma_{\textup{s}}(\fx)$ is a scale factor for the collision operator.
This work focuses on the linear kinetic equation where $Q$ is defined as
\begin{equation}\label{eq_linear_collision_operator}
	Q(f)(\mathbf{v})=\int_{\Stwo} k(\fv_*, \fv) \left[ f(\fv_*)-f(\fv)\right] d\fv_*.
\end{equation}
The collision kernel $k(\fv_*, \fv)$ models the probability of a collision at the point in space $\fx$, given the pre-collision velocity $\fv_*$ and post-collision velocity $\fv$ satisfying $\mathbf{v,v_*}\in \Stwo=\set{\fv\in\mathbb{R}^3: \norm{\fv}=1}$, where $\norm{\cdot}$ denotes the usual Euclidean norm. Finding the kinetic density $f$ that solves
Eq.~\eqref{eq_boltzmann} is a well-posed problem, provided suitable boundary and initial conditions are considered \cite{babovsky1998boltzmann}.
The linear kinetic equation possesses some fundamental structural properties, which are intricately related to the physical processes and its mathematical existence and uniqueness theory~\cite{alldredge2019regularized, pmlr-v162-schotthofer22a, xiao2023relaxnet}. To discuss these properties, we first introduce the notation $\inner{\cdot} = \int_{\Stwo}\cdot \intD \mathbf{v}$ for the integral over the velocity space. 
\begin{definition}[Invariant Range]\label{def_invariant_range}
The kinetic density obeys the physical bounds $f(t,\fx,\mathbf{v})\in B$ for all $t>0$, if $f(0,\fx,\mathbf{v})\in B$ for  $B=[0,\infty)$. In particular, this implies non-negativity of $f$ if the initial condition $f_0$ is non-negative.
\end{definition}
\begin{definition}[Conservation]\label{def_conservation}
There exist quantities $\varphi(\mathbf{v}):\Stwo\rightarrow\mathbb{R}$, known as collision invariants of ${Q}$, such that 
\begin{align}
    \label{eq_collision_invariant}
    \inner{\varphi Q(f)}=0,\quad\forall f\in\text{Dom}( Q),
\end{align}
where $\text{Dom}(Q)\subset B$ denotes the domain of the collision operator $Q$.
The above equation implies the conservation law that takes the form
\begin{align}
    \label{eq_conservation}
    \partial_t \langle \varphi f \rangle + \nabla_{\mathbf x} \cdot \langle \fv \varphi f \rangle = 0.
\end{align} 
\end{definition}

\begin{definition}[Hyperbolicity]\label{def_hyperbolicity}
For each fixed $\fv$, the left-hand side advection operator in Equation (\ref{eq_boltzmann}) is hyperbolic over $(t,\mathbf{x}) \in [0, \infty) \times \mathbf X$.
\end{definition}

\begin{definition}[Entropy Dissipation]\label{def_entropy}
Any twice continuously differentiable, strictly convex function $\eta: D\subseteq \mathbb R \rightarrow \mathbb R$, is an entropy density, such that
{\small
\begin{equation}
  \langle \eta'(g)  Q(g) \rangle \leq 0, \quad \forall g\in \mathrm{Dom}(Q) \ \mathrm{s.t.} \ \mathrm{Range}(g) \subseteq D.
\end{equation}
} 
This inequality yields the local entropy dissipation law
\begin{equation}
	\partial_t\inner{\eta(f)} + \nabla_{\fx}\cdot\inner{\fv\eta(f)}= \inner{\eta'(f)Q(f)} \leq 0,\label{eq_kinetic_entropy_dissipation}
\end{equation}
where $h(f)=\inner{\eta(f)}$ is the entropy of the particle system. 
The system is in equilibrium when the above inequality takes an equal sign, and the corresponding kinetic density is called the equilibrium distribution.
\end{definition}
Due to Property~\ref{def_entropy}, the particle system always strives for a state of minimal mathematical, i.e., maximal physical entropy. The domain $D$ of the entropy density is often consistent with the physical bounds $B$ of the problem.  
In this work, the Maxwell-Boltzmann entropy density $\eta(f)= f\log(f) -f$ with $D=B=[0,\infty)$ is considered.

For a machine learning-based surrogate model of the kinetic equation or one of its operators, preserving the structural properties described above in the numerical simulation is crucial. Especially in multi-physics simulations with coupled PDEs, violation of structural properties may cause instabilities of the corresponding PDE solvers leading to non-physical model-based predictions.

\subsection{Contribution}
This work focuses on constructing a DeepONet-based surrogate model of the collision operator $Q$, see Eq.~\eqref{eq_linear_collision_operator}, called \textit{collision DeepONet}. The collision DeepONet is a modified DeepONet that conserves the collision invariance property of $Q$. Further, the approach conserves the invariance of range and hyperbolicity of the hybrid PDE system obtained by substituting the collision DeepONet into Eq.~\eqref{eq_boltzmann}.
We propose a structured, entropy-based sampling strategy to train the collision DeepONet. This tailored sampling strategy allows to efficiently sample from the domain of $Q$ using the entropy dissipation principle to reject samples that are unlikely to appear in any simulation. A related approach of entropy based sampling has been developed for moment closures~\cite{pmlr-v162-schotthofer22a} and rarefied gas simulations~\cite{XIAO2023112278}. We demonstrate through 
various numerical experiments that the proposed collision DeepONet methods provide robust and interpretable surrogate models.

\subsection{Related work}
Deep Operator Network (DeepONet) techniques have emerged as a promising framework for operator learning. Unlike traditional deep neural networks primarily focusing on learning functions, DeepONets learn operators. 
\citet{lu2019deeponet,lu2021learning} first proposed  DeepONet, building upon the universal approximation theorem for operators \citep{chen1995universal}. It has been applied to various realistic problems, such as fluid flows \citep{cai2021deepm}, bubble growth dynamic prediction \citep{lin2021operator}, and predicting hypersonic shock flows \citep{MR4316010}. Recently, DeepONet and its various modified models \citep{lee2023hyperdeeponet,wang2021learning} have been gaining attention as they show promise to overcome the curse of dimensionality \citep{lanthaler2022error} for high dimensional settings.
The DeepONet structure consists of two (deep) neural networks: the branch-net and the trunk-net. Many studies \citep{hadorn2022shift,lee2023hyperdeeponet, LU2022114778,meuris2023machine} interpret the roles of the branch-net and trunk-net in the DeepONet from the perspective of the basis functions of the target operator's image. 
The trunk-net creates the basis of the target function, while the branch-net generates scales dependent on the function input of the target operator.

Kinetic equations are subject to rising attention in the surrogate modeling community due to their versatility as a mathematical model at the cost of high computational complexity. Physics informed surrogate models for the Boltzmann equation ~\cite{lou2020physicsinformed,li2021physicsinformed,mishra2021physics}, hybrid continuum-kinetic models~\cite{XIAO2023112278,xiao2021using}, and closures of the moment system of kinetic equations~\cite{Han21983, pmlr-v162-schotthofer22a, li_grad_global, huang2021machine, weinan_machine, bois2020neural, Maulik2020neural} is an active focus of research. By the authors' knowledge, provable enforcement of analytical properties of the collision operator has not been pursued yet.

\section{Collision DeepONet}\label{sec_theory}

The high-dimensional joint phase space $(\fx,\fv)$ and the high-dimensional integral of the collision operator $Q$ pose severe challenges for analysis and accurate simulations. Using the deep operator learning approach, we strive to build a surrogate model of the collision operator $Q$ that preserves its core analytical structures.

 To facilitate the discussion, consider the kinetic density $f(\cdot,\cdot,\fv)$ at a fixed time $t$ and space $\fx$ as an integrable function on $\Stwo$, i.e., $f\in L^1(\Stwo, B)$. The collision operator $Q$ maps from $L^1(\Stwo,B)\mapsto L^1(\Stwo)$. The linear collision operator $Q$ considered in this work has a single collision invariant $\varphi=1$ and a single corresponding conserved quantity, namely the mass of the system, denoted by $q(t)=\inner{\varphi f}$. The discussion hereafter will focus on this linear collision operator. However, it is noteworthy that the first proposed method, see~\ref{sec_methodI} below, is straightforwardly extendable to a setting with multiple collision invariants, e.g., for the nonlinear collision operators implying conservation of mass, momentum, and energy.

\subsection{Deep operator learning}
A DeepONet-based surrogate model, in the following denoted by $Q_\theta$,\footnote{The subscript $\theta$ denotes the set of all parameters characterizing the DeepONet approximation of the collision operator.}, approximates the target operator $Q$ by generating a basis of the image space of $Q$, and a coefficient vector to reconstruct the function $Q(f)(\fv)$ from the generated basis \citep{hadorn2022shift,lee2023hyperdeeponet,LU2022114778,meuris2023machine} . The DeepONet therefore consists of two modules: the trunk-net to generate the basis functions and the branch-net to generate the coefficient vector.
The trunk-net and branch-net are constructed as neural networks using a discretization of the velocity space $\Stwo$.

As illustrated in Figure \ref{fig:vanilla_deeponet} in Appendix~\ref{app_theory}, the trunk-net $\tau:\Stwo\mapsto\mathbb{R}^p$ generates $p$ basis functions to span the image of the collision operator $Q$. Each element of the trunk-net $\tau_k(\fv)$, $k=1,\dots, p$, is a velocity-dependent function. The branch-net $\beta:\mathbb{R}^m\mapsto\mathbb{R}^p$ takes discrete values of the kinetic density $f$, called sensor values $\bar{f}=[f(\fv_1),...,f(\fv_m)]\in\mathbb{R}^m$, measured at $m$ sensor points $\fv_1,...,\fv_m\in\Stwo$. The sensor values are mapped onto coefficient values $\beta(\bar{f})=[\beta_1(\bar{f}),...,\beta_p(\bar{f})]\in\mathbb{R}^p$ for the corresponding basis functions given by the trunk-net.

The inner product of the two $p$-length vectors yields the approximated target function $Q_\theta(f)(\fv)$ as
\begin{equation}\label{eq_vanilla_deeponet}
    Q_\theta(f)(\fv)=\sum_{k=1}^p  \beta_k(\bar{f})\tau_k(\fv)+b
\end{equation}
where $b\in\mathbb{R}$ is an additional trainable bias term.

\subsection{Structure preservation}
We reiterate that the collision operator $Q$ determines the conserved quantities in the physical system of interest. Especially when the kinetic model is coupled to other models, physical conservation of, e.g., mass, is expected from the simulated kinetic density. Therefore, the collision DeepONet $Q_\theta$ is also required to preserve the conserved quantities.

In the original DeepONet structure presented in \eqref{eq_vanilla_deeponet}, the learning bias parameter $b$ is a freely trainable value independent of the input $f$, see Fig.~\ref{fig:vanilla_deeponet}. We interpret the trunk-net $\tau$ as the basis of the image of $Q_\theta$. Further, we formally remove the bias term $b$ and introduce an additional trunk-net basis $\tau_{p+1}(\fv)$ with corresponding branch-net coefficient  $\beta_{p+1}(\bar{f})$, see Fig.~\ref{fig:modified_deeponet}.
The corresponding modified DeepONet reads 
\begin{equation}\label{eq_modify_deeponet}
    Q_\theta(f)(\fv)=\sum_{k=1}^{p+1}  \beta_k(\bar{f})\tau_k(\fv).
\end{equation}
In the following, we propose two methods to enforce the preservation of collision invariants of $Q_\theta$, given by Eq.~\eqref{eq_modify_deeponet}, for any input of sensor value $\bar{f}$ by modification of $\tau_{p+1}$ and $\beta_{p+1}$.

\subsubsection{Method \RNum{1}: Orthonogonal trunk-net}\label{sec_methodI}

By Property~\ref{def_conservation}, a collision invariant $\varphi$ is orthogonal to the image of $Q$. Since the trunk-net $\tau$ forms a basis of the image of $Q$, constraining the trunk-net to be orthogonal to  $\varphi$ yields an ansatz to preserve this collision invariant and the corresponding conserved quantities. To this end, we formally interpret the collision invariant $\varphi$ as an additional trunk-net function $\tau_{p+1}=\varphi$ with corresponding fixed branch-net coefficient $\beta_{p+1}=0$. Then, we enforce the orthonormality of $\tau$ via
\begin{align}\label{eq_orthogonality_constr_trunk_net}
    \inner{\tau_{k_1}\tau_{k_2}}&=0,\;\quad\forall 1\leq k_1\neq k_2\leq p+1.
\end{align}
Assuming that Eq.~\eqref{eq_orthogonality_constr_trunk_net} holds true, it follows that
{\small
\begin{align}
\begin{aligned}\label{eq_vanilla_collision_invariant}
      \inner{\varphi Q_\theta(f)}&= \inner{\varphi  \left(\sum_{k=1}^{p+1}  \beta_k(\bar{f})\tau_k\right)}\\
      &=\sum_{k=1}^{p} \beta_k(\bar{f})\inner{\varphi \tau_k} + \beta_{p+1}(\bar{f})\inner{\varphi \tau_{p+1}} = 0.
\end{aligned}
\end{align}
}
The constraint is implemented into the DeepONet training routine by transforming the extended trunk-net $\set{\tau_k}_{k=1}^{p+1}$ into an orthonormal basis $\hat{\tau}\ni\varphi$\footnote{We use the notational convention that $\varphi$ is the basis $p+1$ of $\hat{\tau}$. This can be established by an index switch.} before each forward pass of $Q_\theta$ during training. Naturally, $\hat{\tau}$ fulfills Eq.~\eqref{eq_vanilla_collision_invariant} by construction. Algorithm~\ref{alg_orthonormal_trunk_net} summarizes the training procedure. The training loss here is $\mathcal{L}(Q_\theta,Q)=\frac{1}{T}\sum_{i=1}^T \norm{Q(\bar{f}_i)-Q_\theta(\bar{f}_i)}_2^2$ for a batch of sensor values $\set{\bar{f}_i}_{i=1}^T$ with size $T$.
\begin{algorithm}[t!]

   \caption{Orthonormal trunk-net training (Method \RNum{1})}
   \label{alg:example}
\begin{algorithmic}
   \STATE {\bfseries Input:} DeepONet $Q_\theta$; collision invariant $\varphi$; sampled sensor values $\set{\bar{f}_{i=1}^T}$; sensor points $\fv_i$, $i=1,\dots,m$; number of training epochs $n_{\textup{epoch}}$
   \FOR{$1\leq t\leq n_{\textup{epoch}}$}
      \FOR{$1\leq i\leq T$}
   \STATE Evaluate trunk-net $\tau_k(\fv)$ $\forall k=1,\dots,p$
   \STATE Evaluate branch-net $\beta_k(\bar{f}_i)$ $\forall k=1,\dots,p$
   \STATE  $\set{\hat{\tau}}_{k=1}^{p+1}\gets \textup{Orthonormalize}([\varphi,\tau])$
   \STATE  $Q_\theta(\bar{f}_i)\gets \sum_{k=1}^p  \beta_k(\bar{f}_i)\hat{\tau}_k(\fv)$
    \ENDFOR
   \STATE  Evaluate $  \mathcal{L}(Q_\theta,Q)$
   \STATE     $\theta^{t+1}\gets\textup{Optimization step}(\theta^{t},\nabla_\theta  \mathcal{L}(Q_\theta,Q))$
   \ENDFOR
\end{algorithmic}
\label{alg_orthonormal_trunk_net}
\end{algorithm}

This approach is extendable to collision operators with multiple collision invariants $\varphi$, e.g., $[1,\fv,\norm{\fv}_2^2]$ for a nonlinear collision operator. Since collision invariants are typically pairwise orthonormal, the constraint in Eq.~\eqref{eq_orthogonality_constr_trunk_net} can be directly applied using  Algorithm~\ref{alg_orthonormal_trunk_net}.

We remark, that without changing the internal structure of the DeepONet surrogate $Q_\theta$, a natural approach to including collision invariants is adding a physics-informed loss to the loss function used for the training of $Q_\theta$. The training loss reads then
\begin{align}
\begin{aligned}\label{eq_soft_constraint}
       \mathcal{L}(Q_\theta,Q)=&
       \frac{1}{T}\sum_{i=1}^T \norm{Q(\bar{f}_i)-Q_\theta(\bar{f}_i)}_2^2\\ 
       &+\lambda \sum_{1\leq k_1\neq k_2\leq p+1}\abs{\inner{\tau_{k_1}\tau_{k_1}}}\;,
\end{aligned}
\end{align}
 recalling that the collision invariant for the linear collision operator $\varphi = 1$ is the trunk-net function with index $p+1$. Further, $\lambda>0$ is a regularization parameter. This soft constraint enforces the collision invariance only empirically, without strict guarantees of structure preservation for unseen data.

\begin{algorithm}[t!]
   \caption{Bias adapted trunk-net training (Method \RNum{2})}
\begin{algorithmic}
   \STATE {\bfseries Input:} DeepONet $Q_\theta$; collision invariant $\varphi$; sampled sensor values $\set{\bar{f}_{i=1}^T}$; sensor points $\fv_i$, $i=1,\dots,m$; number of training epochs $n_{\textup{epoch}}$
   \FOR{$1\leq t\leq n_{\textup{epoch}}$}
      \FOR{$1\leq i\leq T$}
   \STATE Evaluate trunk-net $\tau_k(\fv)$ $\forall k=1,\dots,p$
   \STATE Evaluate branch-net $\beta_k(\bar{f}_i)$ $\forall k=1,\dots,p$
   \STATE $\tau_{p+1}(\fv)\gets1$
    \STATE $\beta_{p+1}(\bar{f}_i)\gets -\sum_{k=1}^p  \beta_k(\bar{f}_i) \frac{\inner{\tau_k\varphi}}{\inner{\varphi}}$
   \STATE  $Q_\theta(\bar{f}_i)\gets \sum_{k=1}^{p+1}  \beta_k(\bar{f}_i)\hat{\tau}_k(\fv)$
   \ENDFOR
   \STATE  Evaluate $  \mathcal{L}(Q_\theta,Q)$
   \STATE     $\theta^{t+1}\gets\textup{Optimization step}(\theta^{t},\nabla_\theta  \mathcal{L}(Q_\theta,Q))$
   \ENDFOR
\end{algorithmic}
\label{alg_biased_trunk_net}
\end{algorithm}

\subsubsection{Method \RNum{2}: Bias adaption of the trunk-net}\label{sec_methodII}

It is sufficient to enforce orthonormality only of $\varphi$ with each $\tau_k$ for $k<p+1$, which induces less computational effort than pairwise orthonormality of the entire extended trunk-net. In the following we propose a DeepONet modificatin that incoropartes this idea. 
Consider the collision invariance property corresponding to mass conservation for the modified DeepONet representation~\eqref{eq_modify_deeponet},
\begin{align*}
      \inner{\varphi Q_\theta(f)}= \sum_{k=1}^{p+1}  \beta_k(\bar{f})\inner{\varphi \tau_k}= 0.
\end{align*}
Using the ansatz 
\begin{align*}
   \beta_{p+1}(\bar{f}) =  -\sum_{k=1}^p  \beta_k(\bar{f}) \frac{\inner{\varphi\tau_k}}{\inner{\varphi}}\;.
\end{align*}
for the last branch-net coefficient and $\tau_{p+1}(\fv) =1$ for the last trunk-net basis, it is straightforward to show
\begin{equation*}
\inner{\varphi Q_\theta(f)} =
\sum_{k=1}^{p}  \beta_k(\bar{f})\inner{\varphi \tau_k} -\sum_{k=1}^p  \beta_k(\bar{f}) \frac{\inner{\tau_k\varphi}}{\inner{\varphi}}\inner{\varphi}=0\;,
\end{equation*}
thus enforcing mass conservation by construction. Algorithm~\ref{alg_biased_trunk_net} summarizes the corresponding training procedure.

\begin{algorithm}[t!]

\caption{Entropy-based data sampling - single sample}
\begin{algorithmic}
   \STATE {\bfseries Input:} Entropy threshold $c$; moment basis $\fm$; standard deviation $\sigma$; DeepONet sensor points $\fv_i$, $i=1,\dots,m$
   \REPEAT
   \STATE Sample $\balpha_{\#}\sim\mathcal{N}(\balpha_{\fu_{\textup{eq}},\#},\sigma)$
   \STATE     $\alpha_{1}\gets-\frac{1}{m_1}\big(\log(m_1) +\log( \inner{\exp\left(\balpha_{\#}\cdot{\fm}_\#\right)})\big)$
   \STATE ${f}(\fv_i) \gets \exp(\balpha\cdot \fm(\fv_i))\qquad \forall i=1,\dots,m$
   \UNTIL{$h(f)<c$}
\end{algorithmic}
\label{alg_sampling}
\end{algorithm}

\section{Entropy-based training data sampling}
Data-driven surrogate models for PDEs often use training data generated by high-resolution numerical simulations. Drawbacks of this approach is immense computational cost of gathering data and the sampling bias of the training data driven by the choice of, e.g., simulation geometries, initial and boundary conditions. 

Therefore, we investigate strategies to directly sample the domain $L^1(\Stwo, B)$ of the collision operator $Q$. Sampling naively in an infinite dimensional space has the inherent difficulty that most sampled data is very unlikely to ever appear in a simulation and, therefore, has little value as training data. To quantify the relevancy of training data, consider the entropy dissipation property of the linear kinetic equation.

The entropy dissipation of the particle system, see Property~\ref{def_entropy}, yields an instrument to identify regions of importance to sample training data from, i.e., one should sample kinetic densities in a neighborhood of the equilibrium distribution. Since the entropy $h$ is convex in $f$ with the minimum at the equilibrium distribution, $h(f)$ can be used directly with a threshold parameter $c$ to use the sub-level set 
\begin{align*}
    \mathcal{F}_c=\set{f\in L^1(\Stwo,B)\,|\,h(f)\leq c}
\end{align*}
as the sampling domain for the kinetic density. It remains to determine a suitable sampling ansatz as well as a distribution over $\mathcal{F}_c$.

Possible ways to sample $\mathcal{F}_c$ have been explored in recent works on machine-learning-based surrogate models for kinetic equations and their moment systems, and general transport equations, e.g., sampling Fourier coefficients~\cite{huang2021machine},  deviations from Gaussian distributions~\cite{MR4554720} or the BGK ansatz~\cite{XIAO2023112278}. Particularly relevant for this work is an entropy-closure based sampling approach, that has been explored for the moment system of the linear kinetic equation~\cite{pmlr-v162-schotthofer22a, entropy_AIAA}. To facilitate the discussion, we define the moment vector $\fu=\inner{\fm f}\in\mathbb{R}^n$ of the kinetic density $f$ for a given set of basis functions $\fm=[m_1,\dots,m_n]:\Stwo\rightarrow\mathbb{R}^n$. In this work, the spherical harmonics basis is considered. Specifically, this implies that the conserved quantity mass $q=\inner{m_1 f}$ corresponds to the first moment $u_1$ of the kinetic density. 

The inverse problem of the density-to-moment map $\inner{\fm f}$ can be formulated as the minimal-entropy closure~\cite{alldredge2019regularized},
\begin{align}\label{eq_entropyOCP}
\min_{f\in \mathcal{F}_{\fm}}h(f)\quad  \text{ s.t. } \mathbf u=\inner{\fm f},
\end{align}
where $\mathcal{F}_{\fm}=\set{f\in L^1(\Stwo,B)\,|\,\inner{\fm f}<\infty}$. 
When a minimizer $f_\fu$ exists%
\footnote{Even if $\fu\in\set{\fu\, | \inner{\fm f}=\fu,\, f\in \mathcal{F}_{\fm}}$,  there may not exist a solution to \eqref{eq_entropyOCP} \cite{Junk1999,Hauck2008ConvexDA} in general. In the setting of this work, there always exists a solution.}
for Eq.~\eqref{eq_entropyOCP}, it is unique and takes the form
\begin{align}\label{eq_entropyRecosntruction}
	{f}_\fu(\fv) = \eta'_*(\boldsymbol{\alpha}_\fu\cdot \mathbf m(\fv)),
\end{align}
where the Lagrange multiplier  $ \balpha_{\fu}\in\mathbb{R}^{n} $ is the solution of the convex dual problem
\begin{align}\label{eq_entropyDualOCP}
	\balpha_\fu =  \underset{\balpha\in\mathbb{R}^{n}}{\text{argmin}}
	\,\inner{\eta_*(\balpha\cdot \mathbf m)} - \balpha\cdot \fu 
\end{align}
and $\eta_*$ is the Legendre dual of $\eta$. For $\eta(f)=f\log(f)-f$, the Legendre dual is given by $\eta_*(f)=\exp(f)$.

Equation~\eqref{eq_entropyRecosntruction} provides an ansatz to sample $\mathcal{F}_c$. 
To this end, one samples the Lagrange multiplier $\balpha\sim\mathcal{N}(\mathbf{\mu},{\sigma})$ from a Gaussian distribution with mean $\mathbf{\mu}$ given by the Lagrange multiplier $\balpha_{\fu_{\textup{eq}}}$ corresponding to the equilibrium distribution $f_{\textup{eq}}$ and its moment $\fu_{\textup{eq}}$. The standard deviation $\sigma$ controls how much priority is given to kinetic densities $f$ with high entropy $h(f)$.

Using the collision invariance of $Q$ and the DeepONet surrogate $Q_\theta$, one can simplify the sampling set $\mathcal{F}_c$ by normalization corresponding to the conserved quantities of the kinetic system. 
For the linear kinetic equation, with mass as the conserved quantity, this means that only densities $f$ with unit mass, i.e., $\inner{m_1 f}=1$ are sampled, and the corresponding data-set is denoted by
\begin{align*}
    \overline{\mathcal{F}}_c=\set{f\in L^1(\Stwo,B)\,| \,h(f)\leq c\, \textup{and}\, \inner{m_1 f}=1 }.
\end{align*}
To sample normalized kinetic densities $f$ with $\inner{m_1 f}=1$ for the Maxwell-Boltzmann entropy, consider the ansatz for the first Lagrange multiplier 
\begin{align*}
    \alpha_{\fu,1}=-\frac{1}{m_1}\big(\log(m_1) +\log( \inner{\exp\left(\balpha_{\fu,\#}\cdot{\fm}_\#\right)})\big)
\end{align*}
together with Eq.~\eqref{eq_entropyRecosntruction}, where $\balpha_{\fu,\#}=[\alpha_{\fu,2},\dots,\alpha_{\fu,n}]$ and $\fm_\#(\fv)=[m_{2}(\fv),\dots,m_{n}(\fv)]$. The sampling strategy is summarized in Algorithm~\ref{alg_sampling}.

\section{Inference of the collision DeepONet}
The DeepONet approximation $Q_\theta$ of the linear collision operator is trained on normalized kinetic densities. Thus, at inference time, the input $f$ to $Q_\theta$ needs to be rescaled to unit mass, i.e., $\tilde{f}={f}/{\inner{f}}$.
Since $Q_\theta$ is required to preserve mass (and all other collision invariants), one can simply re-scale the output 
\begin{align}\label{eq_rescale}
    Q_\theta(f) = \inner{f}Q_\theta(\tilde{f})
\end{align}
to obtain the approximation to $Q(f)$. Special consideration is required for the case that $f(\fv)=0$ for all $ \fv\in\Stwo$. The sampling ansatz 
Eq.~\eqref{eq_entropyRecosntruction} can only sample $f>0$, however, $f=0$ is not contained in the training data set. Thus one cannot expect high accuracy of $Q_\theta$ at $f=0$. An important property of the collision operator is $Q(f=0)=0$, i.e. that $f=0$ is an element of the null space of $Q$. To enforce this property, the output of $Q_\theta$ is manually set to zero, if the sensor values $\bar{f}$ are zero-valued.

\begin{table}
\caption{Comparison of the mean relative $L^2$ test error of the proposed conservative DeepONet models for collision operator approximation for $1D$ slab geometry and full velocity space $\Stwo$ with isotropic and Henyey-Greenstein (HG) kernel.}
\label{table:accuracy}
\begin{center}
\resizebox{0.49\textwidth}{!}
{
\begin{tabular}{l|ccccc}
\hline
 &
  Vanilla (w/o bias) &
  Vanilla (w/ bias) &
  \begin{tabular}[c]{@{}c@{}}Soft constraint\end{tabular} &
  \begin{tabular}[c]{@{}c@{}}Orthogonal\\ (Method \RNum{1})\end{tabular} &
  \begin{tabular}[c]{@{}c@{}}Bias adaption\\ (Method \RNum{2})\end{tabular} \\ 
  \hline
1D isotropic & 0.0132 & 0.0146 & 0.0135 & \textbf{0.0084} & 0.0104 \\
1D HG, $g=0.9$ & 0.1064 & 0.1074 & 0.1771 & 0.1056 & \textbf{0.0425} \\
$\Stwo$ isotropic & 0.1086 & 0.1123 & 0.1435 & \textbf{0.0710} & 0.1063 \\
$\Stwo$ HG, $g=0.9$ & 0.2020 & 0.2051 & 0.3234 & \textbf{0.1453} & 0.1950 \\
\hline
\end{tabular}
}
\end{center}
\end{table}

\begin{table}[t!]
\caption{Architectures of the used DeepOnet variants. Displayed are the number of output channels of each fully connected layer of branch- and trunk-nets. Hyperbolic tangent functions are used as the activation functions. 
}
\label{tab_architectures}

\resizebox{0.49\textwidth}{!}
{
\begin{tabular}{c|cc}
\hline
 Velocity space & branch-net  & trunk-net  \\
 \hline
 $1D$ Slab &$100-16-16-16$ & $1-16-16-16$ \\
$\Stwo$  &$800-100-100-16$ & $3-100-100-16$ \\

\hline
\end{tabular}
}
\end{table}

\subsection{Hybrid solver for the kinetic equation} 

In order to generate accurate training data by solving the linear kinetic equation, a standard discretization scheme~\cite{kitrt_steffen} for space, time, and velocity is required that matches the discretization of the collision operator. The velocity space $\Stwo$ is discretized using a tensorized Gauss-Legendre quadrature rule of order $N$, which yields $m$ quadrature points. Applying the discretization to \cref{eq_boltzmann} yields a system of $n_q$ hyperbolic conservation laws with solution vector $f_q(t,\fx)=f(t,\fx,\fv_q)$, coupled by the discretized collision operator $Q_q$, i.e., 
{
\begin{equation}\label{eq_bolzmann_semi} 
	\partial_t f_q+\fv_q \cdot \nabla_{\fx} f_q +\sigma_{\textup{a}} = \sigma_{\textup{s}}Q_q(f_q) + p_q\;,
\end{equation}
}
where $Q_q$ and $p_q$ denote the velocity space discretizations of the corresponding operators for $ 1\leq q\leq m$.
The $n_q$ quadrature points of the velocity space discretization are chosen to coincide with the sensor points $\fv_1,\dots,\fv_m$ of the DeepONet approximation $Q_\theta$, thus the solution vector $f_q$ coincides with the DeepONet sensor values $\bar{f}$.  Exchanging the discretized collision operator $Q_q$ by $Q_\theta$ yields the hybrid system 
\begin{equation}  \label{eq_boltzmann_hybrid}
	\partial_t \bar{f}+\fv_q \cdot \nabla_{\fx} \bar{f} +\sigma_{\textup{a}} = \sigma_{\textup{s}}Q_\theta(\bar{f}) + \bar{p},
\end{equation}
with $\bar{p}=[p_1,\dots,p_m]$. The equation system \cref{eq_boltzmann_hybrid} preserves the hyperbolicity property of \cref{def_hyperbolicity}. If  $Q_\theta$ preserves the collision invariants of $Q$, the hybrid system \cref{eq_boltzmann_hybrid} implies a corresponding conservation law, see \cref{def_conservation}. Invariance of range, see \cref{def_invariant_range} is naturally preserved since $Q_\theta$ in \cref{eq_boltzmann_hybrid} contains a limited number of neurons.
A standard finite volume scheme is employed for the space and time discretization, where space is discretized with a structured grid with $n_\fx$ cells, that resembles the domain of the corresponding test case and the time evolution is computed explicitly by a Runge-Kutta scheme with $n_t$ steps.

\begin{table}
\caption{ Error of the collision invariance property $\abs{\inner{\varphi Q_\theta}}$ with $\varphi=1$ on unseen test data for various DeepONet modifications.}
\label{table:Qdv}
\resizebox{0.49\textwidth}{!}
{
\begin{tabular}{l|ccccccc}
\hline
& Vanilla (w/o bias) & Vanilla (w/ bias) & Soft constraint & \begin{tabular}[c]{@{}c@{}}Orthogonal\\ (Method \RNum{1})\end{tabular} & \begin{tabular}[c]{@{}c@{}}Bias adaption\\ (Method \RNum{2})\end{tabular} \\
 \hline
 1D HG, $g=0.9$ & $1.17\times10^{-3}$ & $1.03\times10^{-3}$ & $1.58\times10^{-3}$ & $4.16\times10^{-7}$ & $\mathbf{7.62\times10^{-8}}$ \\
 $\Stwo$ HG, , $g=0.9$ & $7.58\times10^{-4}$ & $6.75\times10^{-4}$ & $5.38\times10^{-3}$ & $\mathbf{1.51\times10^{-8}}$ & $1.74\times10^{-8}$ \\
\hline
\end{tabular}
}
\end{table}

\section{Numerical experiments}\label{sec_results}
In the following, the performace of the proposed collision DeepONet methods are evaluated through a series of synthetic and simulation experiments. We train all experiments for 10,000 epochs and selected the best model based on the lowest loss achieved during training. We employ the Adam optimizer~\cite{kingma2017adam} with a learning rate of 1e-3 in a full-batch training setting.  

\subsection{Isotropic kernel and Henyey-Greenstein kernel approximation}
The Henyey-Greenstein kernel~\cite{henyey1941diffuse} is an anisotropic model for the scattering kernel $k(\fv_* , \fv)$ of the collision operator in Eq.~\eqref{eq_linear_collision_operator}.
The Henyey-Greenstein kernel reads
\begin{align*}
    k_g(\fv_*,\fv) = \frac{1-g^2}{4\pi(1-2g\fv\cdot\fv_* + g^2)^{3/2}}
\end{align*}
with the anisotropy parameter $g\in[-1,1]$, which allows tuning the scattering from purely backward, i.e., $g=-1$, to isotropic scattering, i.e., $g=0$, and to purely forward scattering $g=1$. 

In the following, we present results for the approximation of the collision operator, with velocity space $\Stwo$ and a simplified one-dimensional setting. 
In general, the collision operator $Q$ is an inherently three-dimensional object, and its one-dimensional counterpart is obtained by assuming a slab-geometry setting with symmetries in the other two spatial directions. The one-dimensional velocity domain is therefore the interval $[-1,1]$.

\begin{figure}
    \centering
    \includegraphics[width=0.49\textwidth]{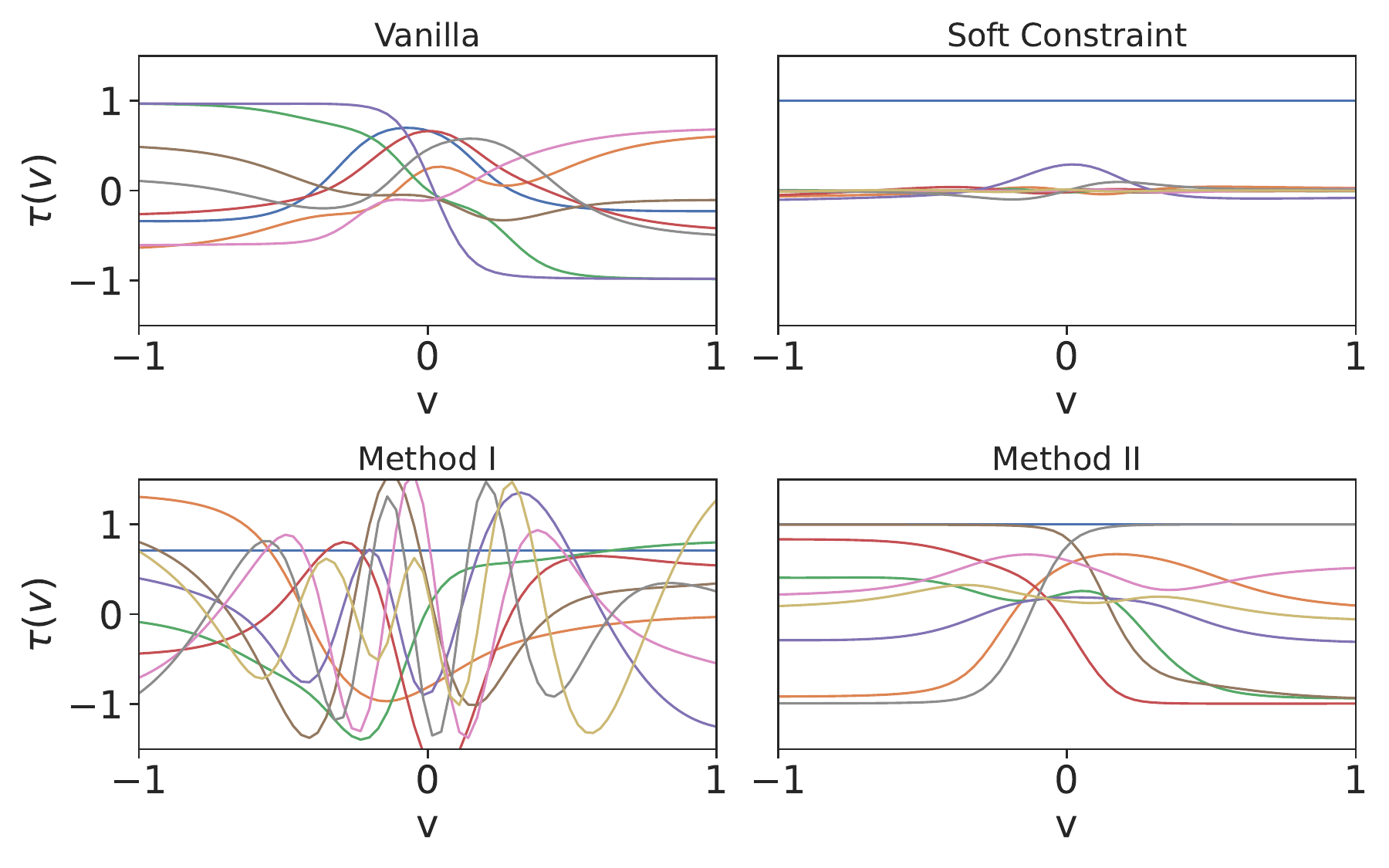}
    \caption{Trunk-net for the 1D Collision operator case with $p=8$ basis functions. The soft orthogonal constraint (upper right) causes the values of $\tau$ to vanish. The orthogonalization with a Gram-Schmidt process (lower left) causes oscillatory behavior. The bias adaption method (lower right)  yields similar trunk-net behavior as the vanilla trunk-net (upper left). }
    \label{fig_truncnet}
\end{figure}

Table \ref{table:accuracy} displays the comparison of the mean relative test errror $\frac{1}{T}\sum_{i=1}^T \norm{Q(\bar{f}_i)-Q_\theta(\bar{f}_i)}_2/\norm{Q(\bar{f}_i)}_2$ of the proposed DeepONet variants, for an isotropic scattering kernel and a Henyey-Greenstein kernel with an-isotropy parameter $g=0.9$ in $1D$ slab geometry and $\Stwo$. 
The DeepONet variants for the $1D$ collision operator have $m=100$ sensor points, and the $\Stwo$ collision operator is equipped with $m=800$ sensor points.
The network architectures --- obtained through an architecture search are displayed in \cref{tab_architectures}. For a given dimension, all DeepONet variants are initialized with the same architecture --- only trunk-net and bias are modified.

All proposed structure-preserving methods train as well as the vanilla DeepONet in the isotropic and forward scattering case. Table \ref{table:Qdv} shows the errors of the collision invariance property $\abs{\inner{\varphi Q_\theta}}$ for 1D and $\Stwo$ velocity space with the HG($g=0.9$) kernel. It indicates that the proposed two methods have smaller errors for the mass conservation property. Figure~\ref{fig_S2_profile} in Appendix~\ref{app_results} exemplary demonstrates the approximation of $Q$ with $Q_\theta$ in the $\Stwo$ velocity space with configuration HG($g=0$).

\begin{figure}
\begin{minipage}{0.24\textwidth}
    \centering
    \includegraphics[width=\textwidth]{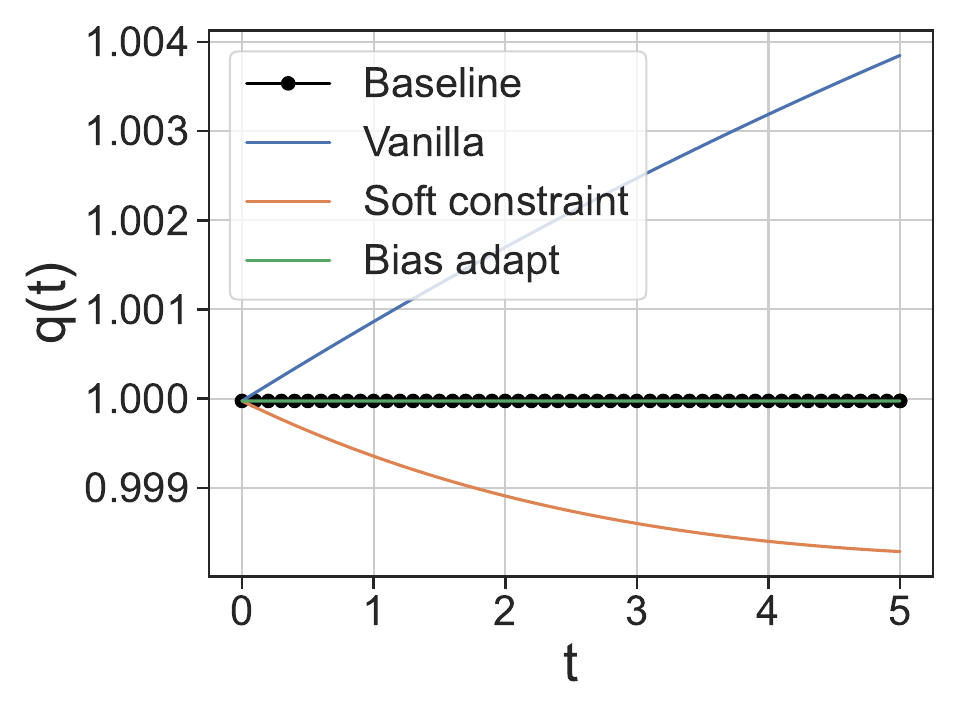}
    \end{minipage}
    \begin{minipage}{0.23\textwidth}
    \centering
          \resizebox{\textwidth}{!}
{
   \begin{tabular}{l|ccc}
    \toprule
    Test case & $n_x$ & $n_t$ & $m$\\ 
    \midrule
    Hom. relaxation & n.a. & $50$ &$100$\\
    Aniso. inflow      & $1000$ & $375$ &$100$\\
    Lattice     & $1600$ & $74$ &$32 $\\
    \bottomrule
        \end{tabular}
}

\end{minipage}

    \caption{Left panel: Homogeneous relaxation test case. The structure-preserving DeepOnet conserves the system's mass, whereas the vanilla versions of DeepONet and the physics-informed version introduce a deviation of the initial mass as time evolves.
Right panel: Discretization settings for the presented test cases with $n_x$ grid cells, $n_t$ time steps, and $m$ DeepONet sensor points.}
    \label{fig_conservation}
\end{figure}

\subsection{Stability of the trunk-net}
The figure~\ref{fig_truncnet} illustrates the different effects of the collision invariant enforcement  on the behavior of the trunk-net. {Enforcing the collision invariance via a physics informed loss, see Eq.~\eqref{eq_soft_constraint}, encourages the DeepONet to scale down the output values of the trunk-net $\tau$, since minimization of $\abs{\inner{\varphi Q_\theta(\bar{f})}}$ is trivial, when $\tau=0 $.} For method I, it is important to use a stable orthogonalization algorithm,  since instable methods, e.g., a Gram-Schmidt process, cause the trunk-net function to oscillate heavily. The bias adaption method has no apparent detrimental influence on the behavior of the trunk-net.

\begin{table}
\caption{An-isotropic inflow test case. Runtime comparison of different collision operator approximation methods.  
}
\label{tab_inflow_timings}
\resizebox{0.49\textwidth}{!}
{
\begin{tabular}{l|cccccc}
\hline
 & Numerical quadrature & Vanilla & Soft constraint & Orthogonal & Bias adaption \\
 \hline
 Run time [iter/s] &1.06 & 41.82 & 39.91 & 2.80 & 40.65 \\
\hline
\end{tabular}
}
\end{table}

\begin{figure}
\centering
\begin{minipage}{0.23\textwidth}
    \includegraphics[width=\textwidth]{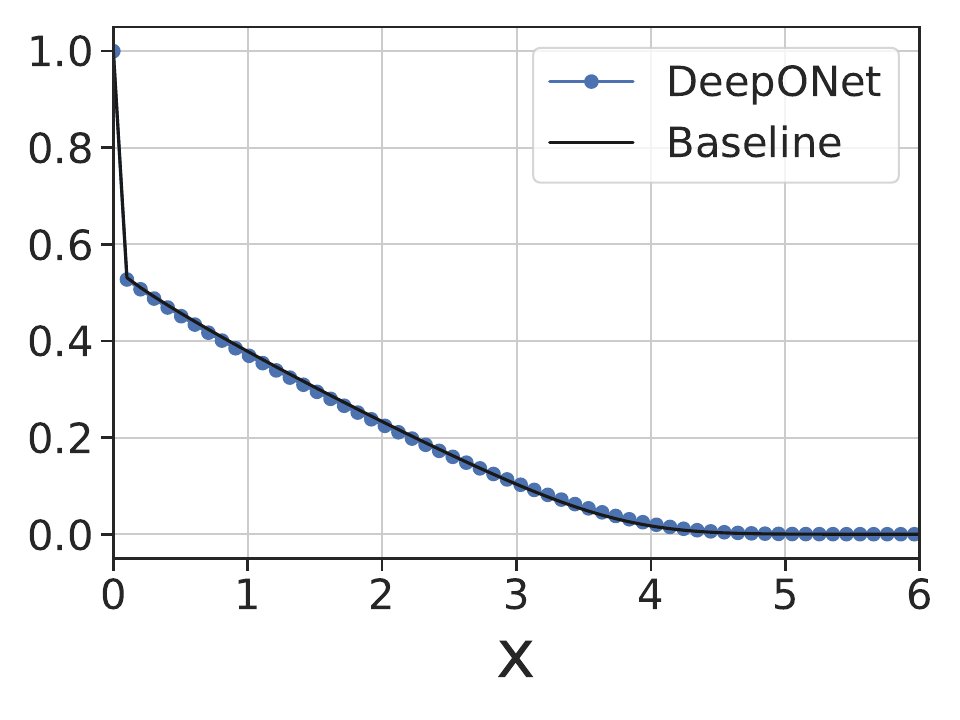}
\end{minipage}
\begin{minipage}{0.23\textwidth}
    \includegraphics[width=\textwidth]{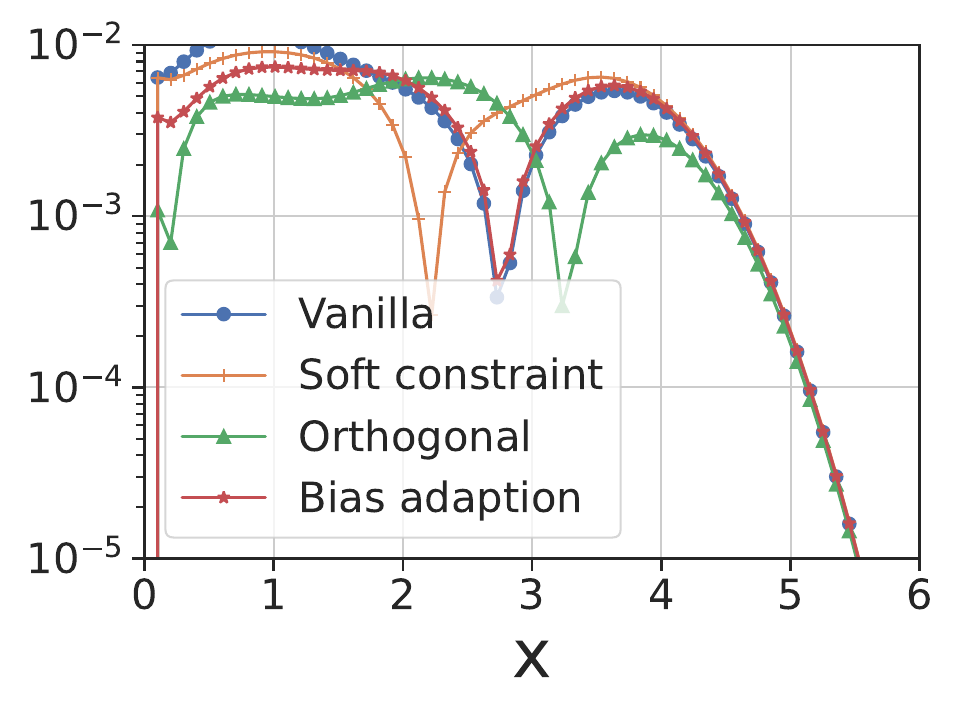}
\end{minipage}
\caption{Anisotropic inflow test case. Left panel: Exact solution and DeepONet approximation through using the bias adaption. Right panel: Comparison of the absolute errors of different DeepONet approximations to the exact solution. 
}
\label{fig_1d_inflow}
\end{figure}

\begin{figure*}
\centering
\begin{minipage}{0.3\textwidth}
        \includegraphics[height=4cm]{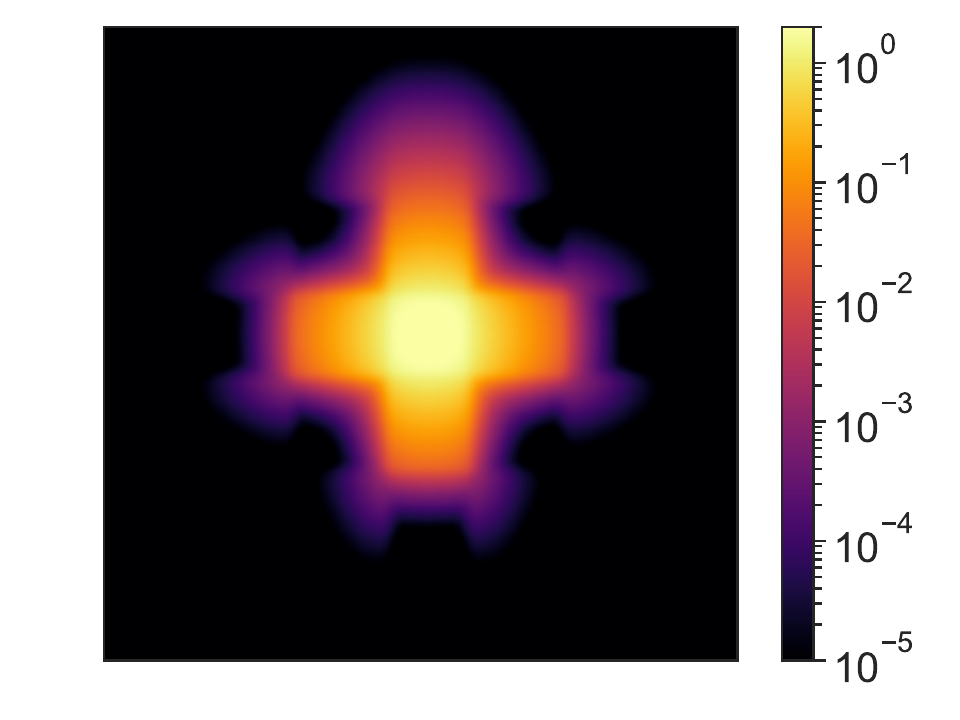}
\end{minipage}
\begin{minipage}{0.3\textwidth}
    \includegraphics[height=4cm]{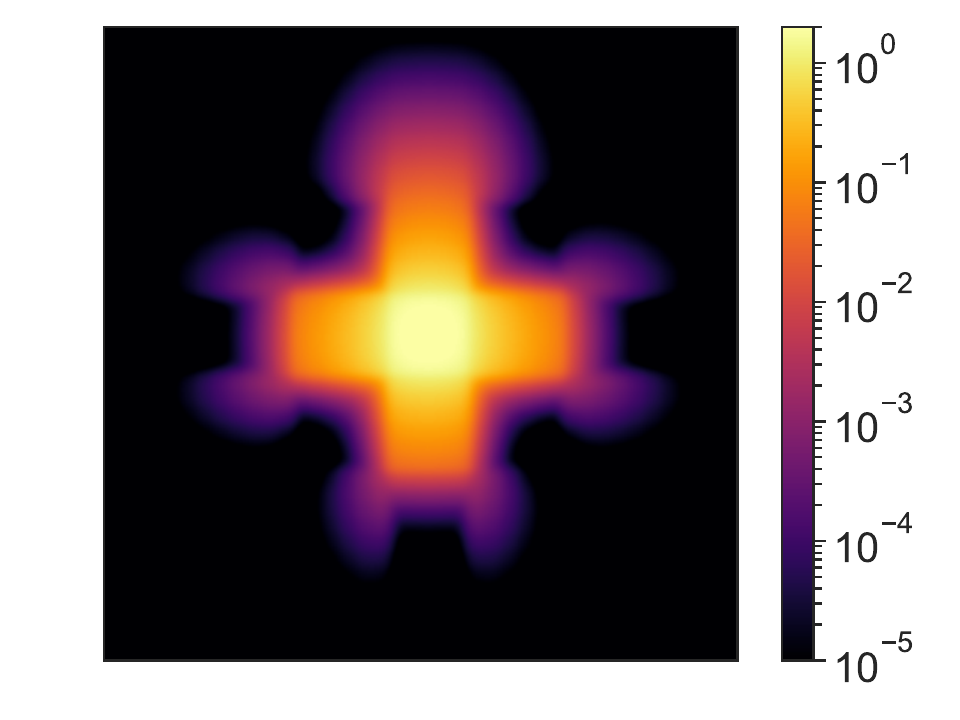}
\end{minipage}
\begin{minipage}{0.3\textwidth}
    \includegraphics[height=4cm]{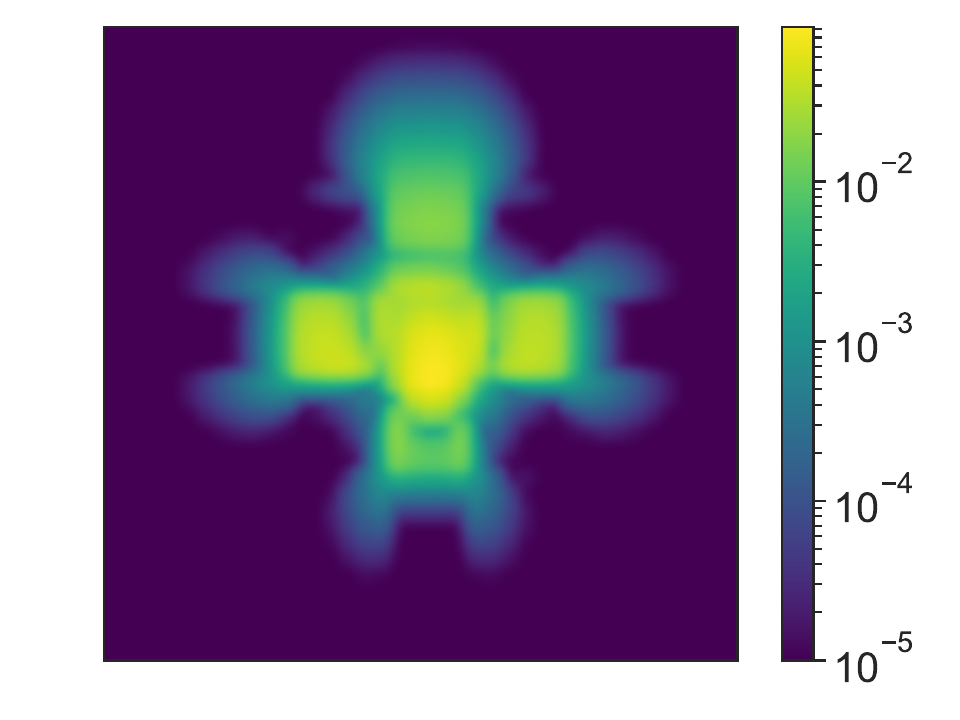}
\end{minipage}

\caption{Lattice test case. Solutions are presented in log scale. Left panel: Baseline solution $q(t_f)$ with $Q_q$ evaluated with a numerical quadrature. Mid panel: Solution $q(t_f)$ with $Q_\theta$ using the bias adaption method. Right panel: Absolute error in each grid cell. The deviation from the baseline method is within reasonable bounds, and critical features such as the blue absorption regions are well resolved.
}
\label{fig_lattice_results}
\end{figure*}

\subsection{Homogeneous relaxation test case}
To verify the conservation property of the collision operator approximation $Q_\theta$, the spatially homogeneous kinetic equation is considered, i.e.
\begin{align}\label{eq_0d_kinetic}
    \partial_t \bar{f} =Q_\theta(\bar{f}),
\end{align}
equipped with initial condition $\bar{f}(t=0,\fv)=\bar{f}_0(\fv)$ for any choice $\bar{f}_0\in L^1(\Stwo,B)$.
The mass, i.e., $q(t)=\sum_{q=1}^m \bar{f}_q$ is the only conserved quantity. Due to the spatial homogeneity and lack of source and absorption terms, no mass is added or subtracted from the system if the collision operator preserves the corresponding collision invariant. \Cref{fig_conservation} displays the time evolution of the mass  $q(t)$ of the system \Cref{eq_0d_kinetic} where $q(t=0)=1$. The proposed structure preserving operator approximation $Q_\theta$ using method II preserves the system's mass $q(t)$, whereas the vanilla DeepONet and the model with physics-informed loss function, see Eq.~\eqref{eq_soft_constraint} introduce a deviation from $q(t=0)$ as time evolves.  This demonstrates the significant outperformance of the proposed modified DeepONet variants in comparison to standard training and soft constraints.

\begin{figure}
\centering
\begin{minipage}{0.23\textwidth}\centering
    \includegraphics[width=0.6\textwidth]{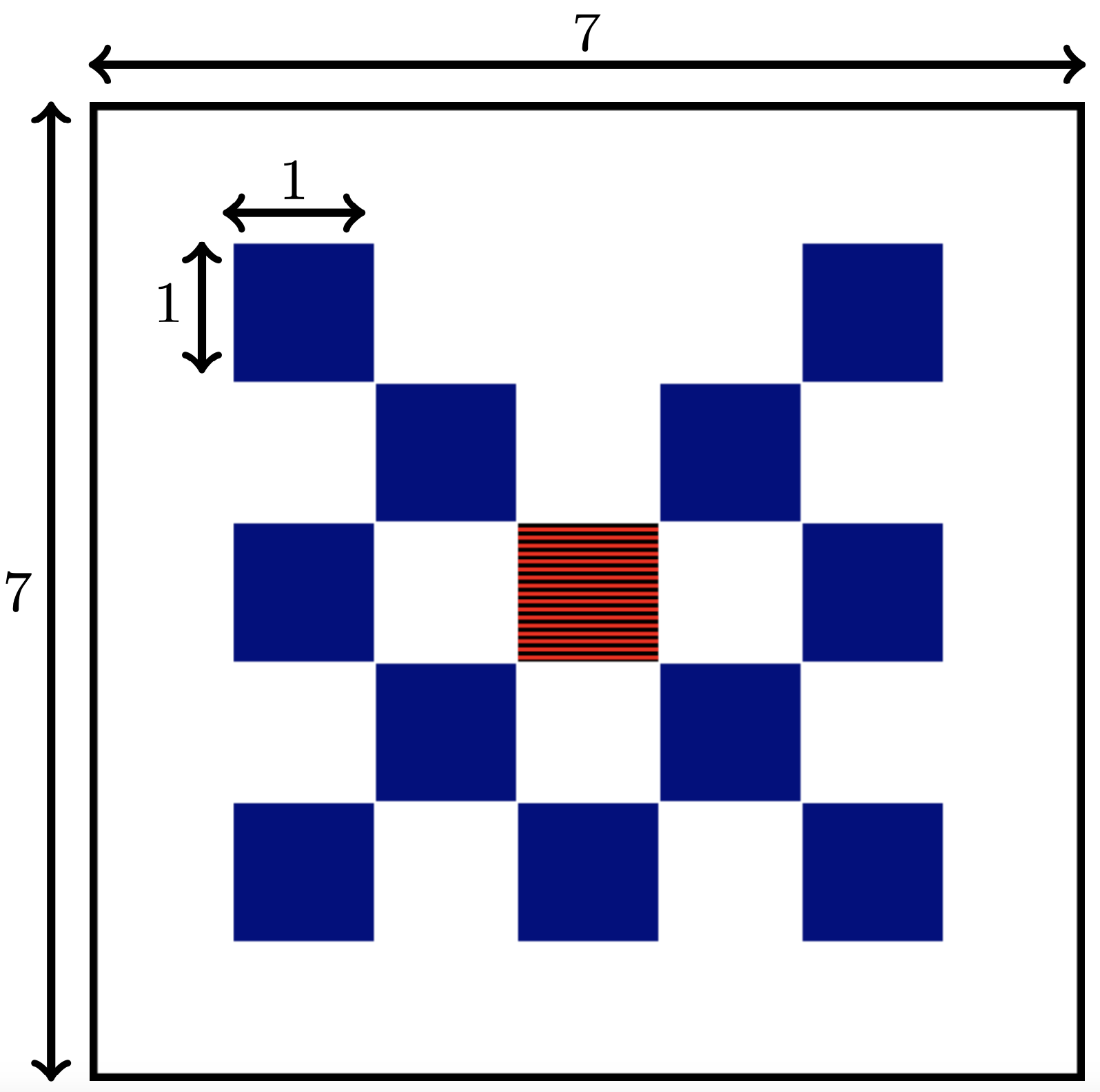}
\end{minipage}
\begin{minipage}{0.23\textwidth}
\label{sample-table}

\vskip 0.15in
\begin{center}
\begin{small}
\begin{sc}
{
\begin{tabular}{l|ccc}
\toprule
Color & $\sigma_{\textup{a}}$ & $\sigma_{\textup{s}}$ & $p$ \\
\midrule
red    & $0$& $1$ & $1$ \\
blue    & $10$& $0$ & $0$\\
white    & $0$& $1$ & $0$ \\
\bottomrule
\end{tabular}
}
\end{sc}
\end{small}
\end{center}
\vskip -0.1in

\end{minipage}

\caption{(Left panel) Geometry of the lattice test case, blue boxes denote high absorption regions and red denotes the radiation source. (Right panel) Numerical values for absorption cross-section $\sigma_{\textup{a}}$, scattering cross-section $\sigma_{\textup{s}}$, and the source term $p$ depending on the test-case geometry. }
\label{fig_lattice_geometry}
\end{figure}

\subsection{An-isotropic inflow test case} 
We test the collision DeepONet variants in the following in a hybrid solver, as defined in Eq.~\eqref{eq_boltzmann_hybrid}.
This simulation test case considers the one-dimensional linear kinetic equation, where the collision operator $Q$  is equipped with a Henyey-Greenstein kernel using $g=0.9$, i.e., a strongly forward scattering setting. No sources or absorption are present, and therefore the setting allows for highly an-isotropic kinetic densities to appear. This corresponds to high entropy values $h(f)$, see ~\cref{eq_entropyOCP}. The data sampling uses an entropy rejection mechanism. Consequently, this test case is tailored to verify the performance of the approach for sensor values $\bar{f}$ that are close to the boundary or even outside the convex hull of the training data set.\footnote{Remark, that the condition $h(f)<c$ defines a convex sublevel set since $h$ is a convex function in $f$.}
The corresponding semi-discrete hybrid equation reads
\begin{align}\label{eq_1d_kinetic}
    \partial_t \bar{f}+ v \partial_{x} \bar{f}
    &=\sigma_{\textup{s}}Q_\theta(\bar{f}),
\end{align}
with $x\in[0,1]$, $v\in [-1,1]$, and $t\in[0,0.7]$. The initial condition is given by $f(t=0,x,v)=0$. 
An an-isotropic inflow condition is imposed at the left boundary of the domain, with $f(t>0,x=0,v) = 0.5$ 
      if $v>0$ and else set to zero,
and the right-hand-side boundary is equipped with a Dirichlet boundary condition $f(t>0,x=1,v)=0$.
The spatial domain is discretized with a grid of size $\Delta x=0.01$ and a corresponding time discretization using CFL number $0.4$. \cref{fig_1d_inflow} demonstrates that all proposed methods yield a reasonable approximation error.
The simulation time of the various DeepONet approximations is compared in \Cref{tab_inflow_timings}, which displays the average iterations per second of the numerical solver. As demonstrated, proposed method I accelerates the standard simulation by a factor $2.2$, and method II has a speedup factor of almost $40$, on par with a non-modified DeepONet. This demonstrates, that the proposed methods are both highly efficient while preserving conserved quantities.

\subsection{Lattice test case}
The Lattice test case~\cite{brunner2005two} is a well-known benchmark for radiation transport simulation. The test case geometry, see Fig.~\ref{fig_lattice_geometry} mimics a nuclear reactor block with a strong radiative source in the domain center and several highly absorptive regions placed in a checkerboard pattern around it.
The corresponding time-dependent hybrid equation is given by ~\Cref{eq_boltzmann_hybrid}
for  $ \fx\in  \mathbf{X} = [0,7]^2$, $t\in[0,3)$.
The equation is equipped with Dirichlet boundary conditions 
$f(t,\fx,\fv)=0$ for $\fx\in\partial \mathbf{X}$ and initial condition  $f(t=0,\fx,\fv)=0$ to obtain a well-posed problem. Furthermore, the scattering kernel $k$ and source term $p$ are assumed to be isotropic and constant in time. The space-dependent absorption cross-section $\sigma_{\textup{a}}$, scattering cross-section $\sigma_{\textup{s}}$, and source term $p$ are given by Fig.~\ref{fig_lattice_geometry}.
The conserved qantity $q=\inner{f}$ at the final time $t_f$ is plotted over the computational grid in \cref{fig_lattice_results}. The geometric features of the source and absorption regions of the lattice test case are captured well by the hybrid method, as the comparison with the baseline demonstrates. The steep drops in radiation density from the red source region to the inner absorption regions is a challenging scenario, yet well resolved by the hybrid method.
\section{Summary}
This work presents the \textit{collision DeepONet} methods, a set of structure preserving DeepONet-based surrogate models for the collision operator of the Boltzmann equation. Both methods are based upon interpretation of the trunk-net as the basis of the image under the DeepONet operator; they use enforcement of orthogonality to enable preservation of collision invariants. Further, an entropy based data-sampler is proposed. Remark, that the current methods are tailored to the velocity space $\Stwo$ and not directly applicable for the nonlinear Boltzmann equation with unbounded velocity domain. This limitation and preservation of the entropy dissipation property is subject to future work, especially with application in rarefied gas dynamics in mind.

\section{Impact statement}
This paper presents work whose main goal is the creation of a structure-preserving deep operator network-based surrogate model for the collision operator of the linear kinetic equation. As in the majority of deep learning research, there are potential societal consequences of our work, none of which we feel must be specifically highlighted here. Regarding ethical aspects, we feel nothing has to be added.


\section*{Acknowledgements}
Jae Yong Lee was supported by a KIAS Individual Grant (AP086901) via the Center for AI and Natural Sciences at Korea Institute for Advanced Study and by the Center for Advanced Computation at Korea Institute for Advanced Study.

The work of Steffen Schotthöfer is sponsored by the Office of Advanced Scientific Computing Research, U.S. Department of Energy, and performed at the Oak Ridge National Laboratory, which is managed by UT-Battelle, LLC under Contract No. DE-AC05-00OR22725 with the U.S. Department of Energy. The United States Government retains and the publisher, by accepting the article for publication, acknowledges that the United States Government retains a non-exclusive, paid-up, irrevocable, world-wide license to publish or reproduce the published form of this manuscript, or allow others to do so, for United States Government purposes. The Department of Energy will provide public access to these results of federally sponsored research in accordance with the DOE Public Access Plan (http://energy.gov/downloads/doe-public-access-plan).

The work of Steffen Schotthöfer and Martin Frank is funded by the Priority Programme ``Theoretical Foundations of Deep Learning (SPP2298)''  by the Deutsche Forschungsgemeinschaft. 


\bibliography{literature}
\bibliographystyle{icml2024}

\newpage
\appendix
\onecolumn

\section{Notations}
The list of notations used throughout the paper is provided in Table \ref{table:notations}.

\begin{table}[h]

\caption{Notations}
\label{table:notations}
\begin{center}
\begin{tabular}{l|l}
\toprule
{ Notation}
&Description\\
\midrule 
$f(t,\fx,\fv)\in L^1(\Stwo,B)$ & The kinetic density in physical space  $\fx$ and velocity space $\fv$ at time $t$ \\
$\fx\in\mathbf{X}$ & Spatial variable  \\
$\fv\in\Stwo$ & Velocity variable \\
$t\geq 0$ & Time variable \\
$\mathbf{X}\subset\mathbb{R}^d$ &$d$-dimensional physical space (we have $d=1,2$ in the current work) \\
$\Stwo=\set{\fv\in\mathbb{R}^3: \norm{\fv}=1}$ & Velocity sphere \\
$Q(f):L^1(\Stwo,B)\mapsto L^1(\Stwo)$ & Collision operator \\
$Q_\theta(f):L^1(\Stwo,B)\mapsto L^1(\Stwo)$ & DeepONet-based surrogate model of the collision operator \\
$k(\fv_*, \fv)$ & Collision kernel \\
$k_g(\fv_*, \fv)$ & Henyey-Greenstein collision kernel with anisotropy parameter $g$ \\
$B=[0,\infty)$ & Physical bound of the kinetic equation \\
$\varphi:L^1(\Stwo)$ & Collision invariants of $Q$ \\
$\eta:D\rightarrow\mathbb{R}$ & Entropy density with domain $D=B$ \\
$h:L^1(\Stwo,B)\rightarrow\mathbb{R}$ & Kinetic entropy \\
$\tau:\Stwo\mapsto\mathbb{R}^p$ & Trunk-net of the DeepONet \\
$\beta(\bar{f})=[\beta_1(\bar{f}),...,\beta_p(\bar{f})]\in\mathbb{R}^p$ & Branch-net of the DeepONet \\
$\theta$ & Set of all trainable paramters of a DeepONet model \\
$\bar{f}=[f(\fv_1),...,f(\fv_m)]\in\mathbb{R}^m$ & Sensor values of the DeepONet \\
$\fv_1,...,\fv_m\in\Stwo$ & Sensor points of the DeepONet.\\
$\mathcal{L}(Q_\theta,Q)$ & Loss function of the training process of a DeepONet \\
$T$ & Training-set size \\
$\fm=[m_1,\dots,m_n]:\Stwo\rightarrow\mathbb{R}^n$ & Basis functions for the entropy closure \\
$\fu=\inner{\fm f}\in\mathbb{R}^n$ & Moment vector of a kinetic density \\
$q=\inner{m_1 f}$ & Mass (conserved quantity) in the linear kinetic system \\
$\balpha_\fu\in\mathbb{R}^n $ & Lagrange multiplier of the minimal entropy closure \\
$\sigma_\textup{a}(\fx)$ & Absorption coefficient in the kinetic equation \\
$\sigma_\textup{s}(\fx)$ & Scattering coefficient in the kinetic equation \\
$p(\fx,\fv)$ & Source term in the kinetic equation \\
\bottomrule
\end{tabular}
\end{center}
\end{table}

\section{Additional material for Section~\ref{sec_theory}}\label{app_theory}
\begin{figure}[t]
\centering
 \begin{minipage}{0.49\textwidth}	\centering
\includegraphics[width=\textwidth]{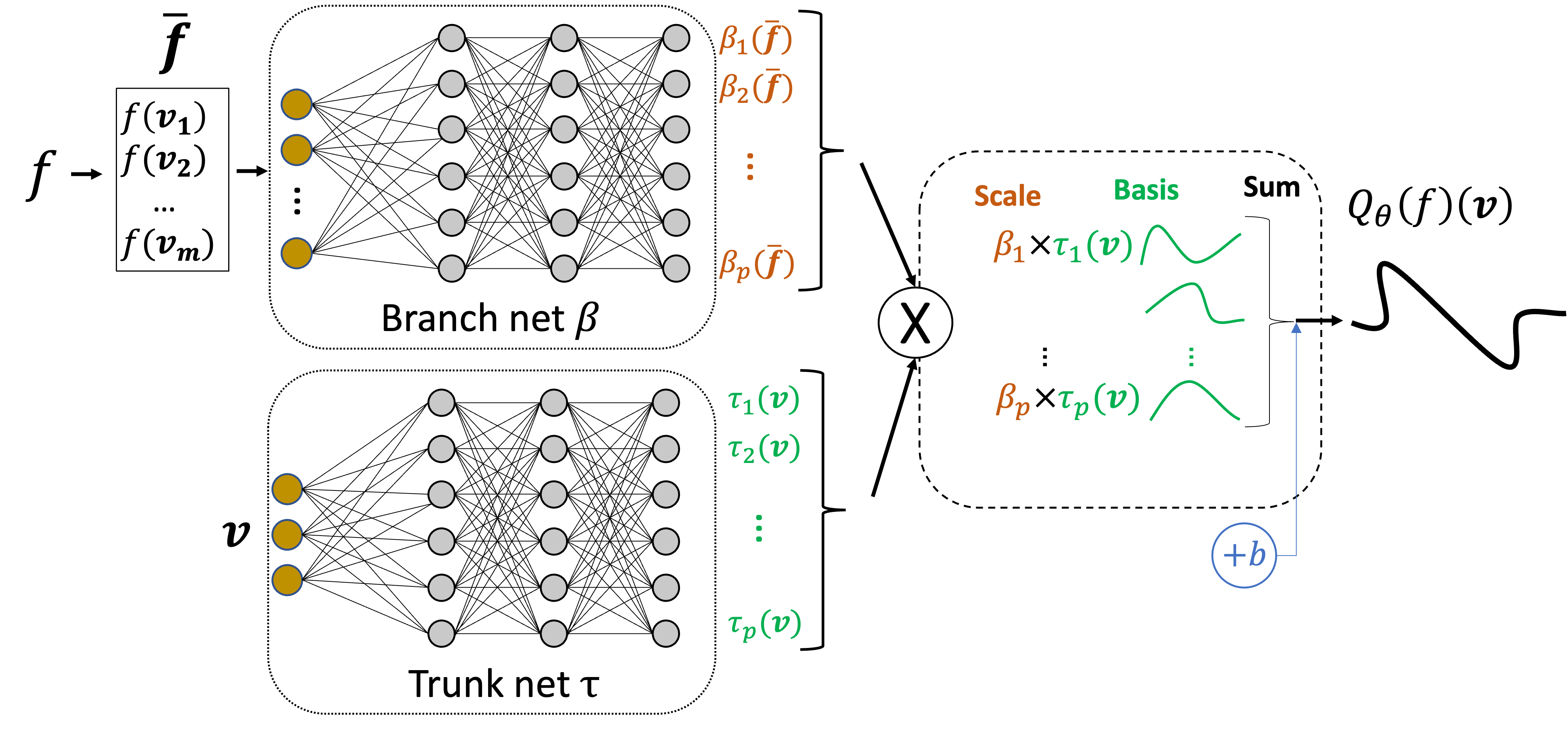}
\caption{The structure of vanilla DeepONet \eqref{eq_vanilla_deeponet}. This original model is unstacked DeepONet with an additional bias proposed in \citet{lu2021learning}.}
\label{fig:vanilla_deeponet}
\end{minipage}
\hfill
 \begin{minipage}{0.49\textwidth}	\centering
\includegraphics[width=\textwidth]{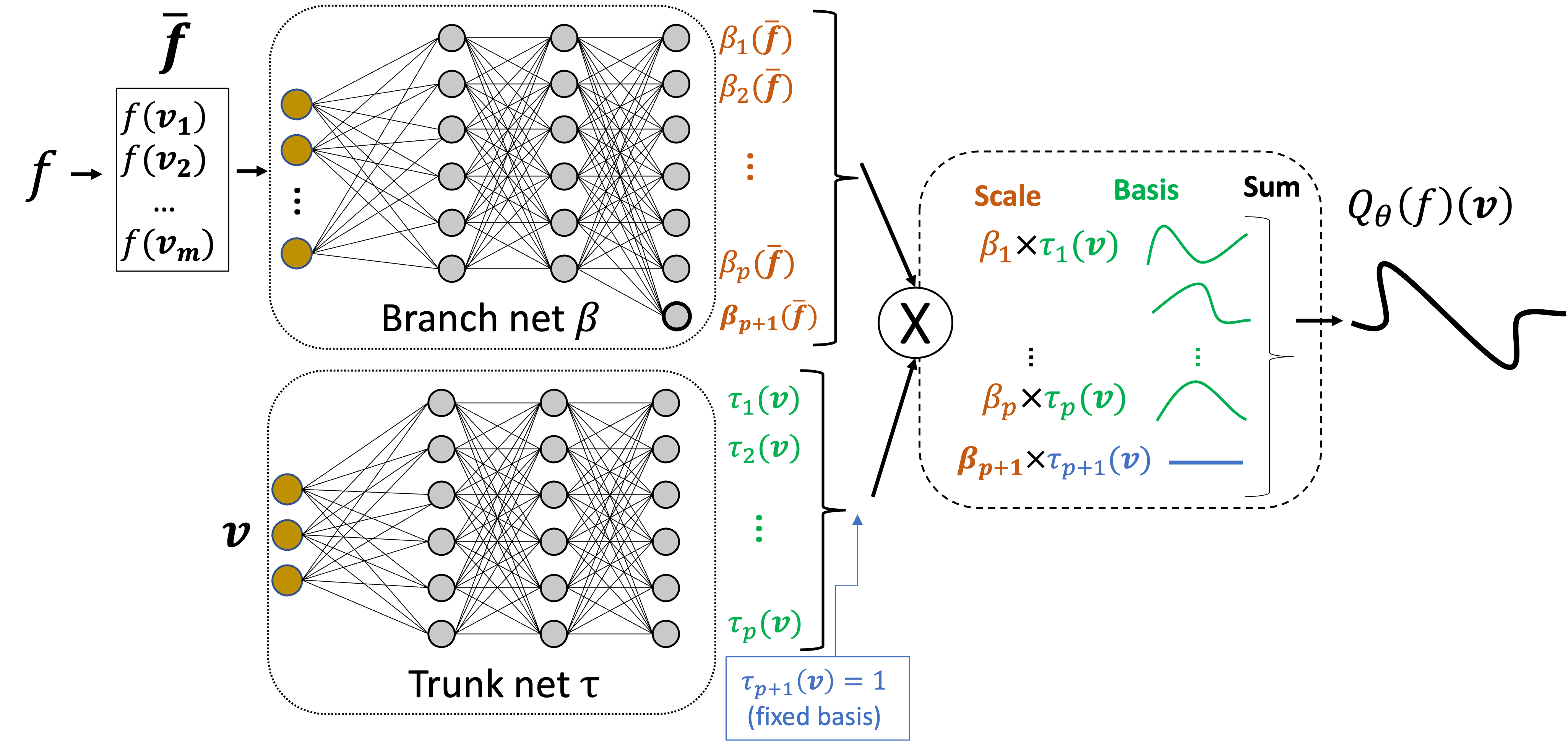}\caption{The structure of modified DeepONet \eqref{eq_modify_deeponet} where the trunk-net is extended by the collision invariant $\varphi=1$.}\label{fig:modified_deeponet}
\end{minipage}
\end{figure}

Figure~\ref{fig:vanilla_deeponet} and Figure~\ref{fig:modified_deeponet} demonstrate the DeepONet architecture and the modification to the trunk-net and bias term. The trunk-net is extended by the number of colission invariants. The bias term of the DeepONet, interpretable as a constant trunk-net function with coefficient $1$, is removed. This allows the description of a orthogonal trunk-net, see Method I in Section~\ref{sec_methodI} that spans the image under the DeepONet operator.  In Method II, the trunk-net is enforced to be only orthorgonal to the collision invariant $\phi$ itself, see Section~\ref{sec_methodII}.

\section{Additional material for Section~\ref{sec_results}}\label{app_results}
\begin{figure}
\centering
{
\includegraphics[width=15cm]{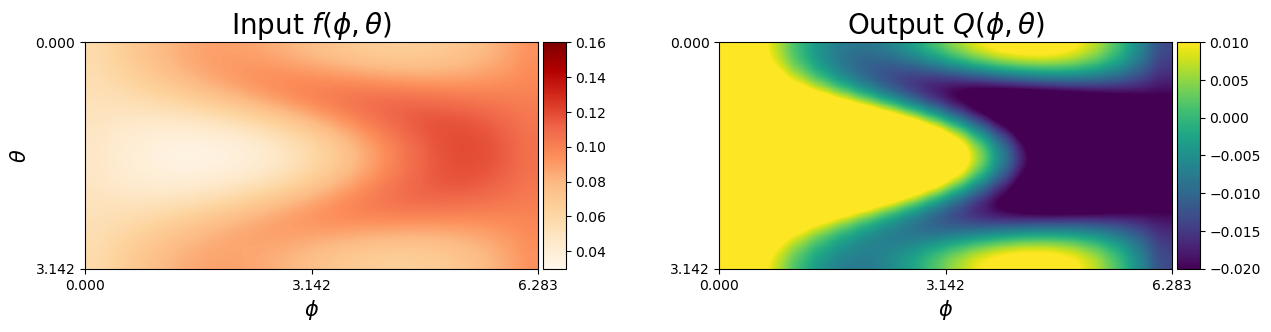}
}
\\
{
\includegraphics[width=15cm]{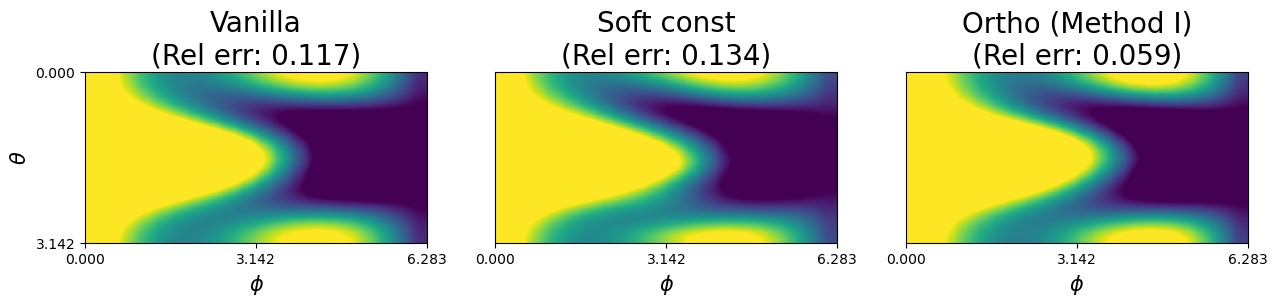}
}
\caption{The approximation of $Q$ with $Q_\theta$ in $\Stwo$ HG case with $g=0$. The domain $\Stwo$ is resolved by $\fv=(\phi,\theta)$. The dynamics of the operator are well resolved by the orthogonal DeepONet, whose relative error is significantly smaller than the error of the baseline and soft constraint versions.}
\label{fig_S2_profile}
\end{figure}

In Figure~\ref{fig_S2_profile}, the effect of the linear collision operator with isotropic collision kernel is demonstrated and compared to various DeepONet approximations. The velocity space $\Stwo$ is parameterized in spherical coordinates with variables $ 0\leq\theta\leq\pi$ and $\leq\phi<2\pi$. 

\end{document}